\documentclass[amsfonts]{article}
\usepackage{amssymb}

\newtheorem{theorem}{Theorem}[section]
\newtheorem{corollary}[theorem]{Corollary}

\newtheorem{proposition}[theorem]{Proposition}
\newtheorem{lemma}[theorem]{Lemma}
\newtheorem{definition}[theorem]{Definition}
\newtheorem{remark}[theorem]{Remark}

\def\bR{\mathbb{R}}
\def\bZ{\mathbb{Z}}

\def\cA{\mathcal{A}}
\def\cB{\mathcal{B}}

\def\cD{\mathcal{D}}

\def\cF{\mathcal{F}}

\def\cK{\mathcal{K}}
\def\cP{\mathcal{P}}
\def\cM{\mathcal{M}}

\begin{document}

\title{{\bf SPECIAL ISSUE - LIMIT
THEOREMS AND TIME SERIES} \\A Cluster Limit Theorem for Infinitely Divisible Point Processes}

\author{Raluca Balan\footnote{ Department of Mathematics and Statistics, University of Ottawa,
585 King Edward Avenue, Ottawa, ON, K1N 6N5, Canada. E-mail address:
rbalan@uottawa.ca}\ \footnote{Supported by a grant from the Natural
Sciences and Engineering Research Council of Canada.} \ and Sana
Louhichi\footnote{Laboratoire Jean Kuntzmann, Tour IRMA IMAG.
51 rue des Mathématiques
BP. 53 F-38041 Grenoble Cedex 9. France. E-mail address:
Sana.Louhichi@imag.fr}\, {\footnote{Partially supported by the ANR grant
ANR-08-BLAN-0314-02.}}}
\date{}
\maketitle
\vspace{-0.5cm}
\begin{abstract}
\noindent In this article, we consider a sequence $(N_n)_{n \geq 1}$
of point processes, whose points lie in a subset $E$ of $\bR
\verb2\2 \{0\}$, and satisfy an asymptotic independence condition.
Our main result gives some necessary and sufficient conditions for
the convergence in distribution of $(N_n)_{n \geq 1}$ to an
infinitely divisible point process $N$. As applications, we discuss
the exceedance processes and point processes based on regularly
varying sequences.
\end{abstract}

%60G55=point processes
%60G70=extreme value theory; extremal processes

\section{Introduction, notation and main assumptions}

Let $E$ be a locally compact Hausdorff space with a countable basis
(abbreviated LCCB) and $M_p(E)$ be the set of Radon measures on $E$
with values in $\bZ_+$. (Recall that a measure $\mu$ is Radon if
$\mu(B)<\infty$ for any $B \in \cB$, where $\cB$ is the class of
relatively compact Borel sets in $E$.) The space $M_p(E)$ is endowed
with the topology of vague convergence. (Recall that $\mu_n
\stackrel{v}{\to} \mu$ if $\mu_n(f) \to \mu(f)$ for any $f \in
C_K^{+}(E)$, where $\mu(f)=\int_{E}fd\mu$ and $C_K^{+}(E)$ is the
class of continuous functions $f:E \to [0,\infty)$ with compact
support.)

A measurable map $N:\Omega \to M_p(E)$ defined on a probability
space $(\Omega,\cK,P)$ is called a point process. The law of $N$ is
uniquely determined by its Laplace functional:
$L_{N}(f)=E(e^{-N(f)}),f \in \cF$, where $\cF$ is the class of
measurable functions $f:E\to [0,\infty)$ (see e.g. \cite{kall83}).
%\cite{daley-verejones03}, \cite{daley-verejones08}).

A point process $N$ is infinitely divisible (ID) if for any integer
$n \geq 1$, there exist some i.i.d. point processes
$(N_{i,n})_{1\leq i\leq n}$ such that
$N\stackrel{d}{=}\sum_{i=1}^{n}N_{i,n}$. An ID point process enjoys
many properties similar to those of ID random vectors. In
particular, a point process $N$ is ID if and only if there exists a
(unique) measure $\lambda$ on $M_p(E)$ (called the canonical
measure), satisfying $\lambda(\{o\})=0$, ($o$ being the null measure
on $E$), and
\begin{equation}
\label{cond-lambda} \int_{M_p(E)}
(1-e^{-\mu(B)})\lambda(d\mu)<\infty \quad \forall B \in \cB,
\end{equation}
such that:
\begin{equation}
\label{Laplace-N-1} L_{N}(f)=\exp\left\{-\int_{
M_p(E)}(1-e^{-\mu(f)})\lambda(d\mu) \right\}, \quad \forall f \in
\cF.
\end{equation}

A sequence $(N_n)_{n \geq 1}$ of point processes converges in
distribution to $N$ if their laws converge weakly (in $M_p(E)$) to
the law of $N$. We write $N_n \stackrel{d}{\to}N$.

The study of point process convergence is important from the
theoretical point of view (by applying the continuous mapping
theorem, one can obtain various limit theorems for some functionals
of the points see for instance \cite{tyran09} and \cite{BKS10}) and practical (see for instance \cite{leadbetter91}, \cite{davis-mikosch98}, \cite{MJS07} and the
references therein). Point processes of exceedances play an important role in extreme value theory and their limiting behavior has been extensively studied (we refer for instance to \cite{BNX02}, \cite{bellanger-perera03}, \cite{HHL88}, \cite{husler-schmidt96}, \cite{leadbetter-rootzen88} and the references therein).

The purpose of the present article is to establish {\it minimal}
conditions for the
 convergence $N_n \stackrel{d}{\to}N$, when $N$ is an ID point
process and $E$ is a subspace of $\bR \verb2\2 \{0\}$. This question
has been studied by
several authors in different contexts (see e.g. \cite{DH95}, \cite{HHL88}).
Our contribution consists in providing a general (unifying) result,
which contains the results of \cite{DH95} and \cite{HHL88},
established in two different situations.

The following definition introduces an asymptotic independence condition.

\begin{definition}
{\rm A sequence $(N_n)_{n \geq 1}$ of point processes satisfies {\bf
condition (AI)} if there exists a sequence $(N_{i,n})_{1\leq i\leq k_n, n \geq 1}$ of point processes such that $k_n \to \infty$, 
\begin{equation}
\label{null-array} \lim_{n \to \infty}\max_{i \leq
k_n}P(N_{i,n}(B)>0)=0, \quad \forall B \in \cB,
\end{equation}
and
$$\lim_{n \to \infty}|E(e^{-N_n(f)})-
\prod_{i=1}^{k_n}E(e^{-N_{i,n}(f)})|=0, \quad \forall f \in
C_{K}^{+}(E).$$ }
\end{definition}

Condition (AI) requires that $N_n$ behaves asymptotically (in
distribution) as the superposition of $k_n$ independent point
processes. This condition contains various conditions
encountered in the literature related to asymptotic results for
triangular arrays  (see condition $\Delta(\{u_n\})$ of \cite{HHL88},
condition $\cA(\{a_n\})$ of \cite{DH95}, or condition
(AD-1) of \cite{BL08}).
More precisely, if $N_n=\sum_{j=1}^{n}\delta_{X_{j,n}}$, where
 $(X_{j,n})_{1 \leq j\leq n, n \geq 1}$ is a triangular array of
 random variables, then (AI) holds under suitable dependence conditions imposed
 on the array (see e.g. \cite{BL08}).

In the present article, we relax the independence assumption between
the components $( N_{i,n})_{i \leq k_n}$, by requiring that (AI)
holds. The following result is an immediate consequence of Theorem
6.1 of \cite{kall83}.

\begin{theorem}\label{kall-th} Let $(N_n)_{n \geq 1}$ be a sequence of point
processes satisfying (AI).
 Then $N_n \stackrel{d}{\to}$ some $N$ if
and only if there exists a measure $\lambda$ on $M_p(E)$ satisfying
$\lambda(\{o\})=0$ and (\ref{cond-lambda}), such that
\begin{equation}
\label{conv-int-lambda-n} \int_{M_p(E)}
(1-e^{-\mu(f)})\lambda_n(d\mu) \to \int_{M_p(E)}
(1-e^{-\mu(f)})\lambda(d\mu), \quad \forall f \in C_{K}^{+}(E),
\end{equation}
where $\lambda_n=\sum_{i=1}^{k_n}P \circ N_{i,n}^{-1}$. In this
case, $N$ is an infinitely divisible point process with canonical
measure $\lambda$.
\end{theorem}

This article is organized as follows. Section 2 is dedicated to our
main result (Theorem \ref{DH-th}) and its consequences for
triangular arrays of random variables. This result is based on
Theorem \ref{kall-th} and gives some necessary and sufficient
conditions for the convergence of $(N_n)_{n \geq 1}$ to an ID point
process $N$, in terms of the asymptotic behavior of the pair
$\{(Y_{i,n},N_{i,n})\}_{1 \leq i \leq k_n}$, where $(N_{i,n})_{1
\leq i \leq k_n}$ are given by (AI), and $Y_{i,n}$ is the largest
point of $N_{i,n}$ (in modulus). We apply our main result to triangular arrays of the form $(\xi_i/a_n)_{i,n}$ (see Proposition \ref{DH-th-cor} below). A relationship
between the canonical measure $\lambda$ of $N$ and the extremal index of the sequence $(|\xi_i|)_i$ when it exists, is established in Corollary \ref{extremal-index-cor}.
In Section 3, we apply our result
to the case of exceendace processes and processes based on regularly
varying sequences, and thus recover the results of \cite{HHL88} and
\cite{DH95}, respectively.

\section{The result}

In the present article, we assume that $E$ is one of the following
spaces:\\
(i) $E=[-b,-a) \cup (a,b]$, $E=(a,b]$ or $E=[-b,-a)$ for some $0 \leq a <b<\infty$; \\
(ii) $E=[-\infty,-a) \cup (a,\infty]$, $E=(a,\infty]$ or
$E=[-\infty,-a)$ for some $0 \leq a <\infty$.

Note that in both cases, $E$ is a LCCB subspace of
$\overline{\bR}=\bR \cup \{\pm \infty\}$, which satisfies the
following conditions:
\begin{equation}
\label{compactness-cond} \mbox{for any $x>0$,} \ [-x,x]^c: =\{y \in
E; |y|>x\} \ \mbox{is relatively compact in} \ E,
\end{equation}
\begin{equation}
\label{compactness-cond2} \mbox{for any compact set $K\subset E$,
there exists $x>0$ such that $K \subset [-x,x]^c$.}
\end{equation}

\vspace{3mm}

\noindent Let $N$ be an infinitely divisible point process on $E$,
with canonical measure $\lambda$. We assume that:
\begin{equation}
\label{cond-N-pm-infty} \mbox{if the space $E$ contains $\pm
\infty$, then $N(\{\pm \infty\})=0$ a.s.}
\end{equation}

\begin{remark}
{\rm (a) The Poisson process $N$ on $[0,1]$ of intensity $1$ is
included in our framework. In this case, we can exclude $0$ from the
space, since the points of the process are strictly positive.
$N$ is a well-defined point process on $E=(0,1]$.

(b) The Poisson process on $[0,\infty)$ of intensity $1$ is {\em
not} included in our framework. It is known that $N=\sum_{\i \geq
1}\delta_{\Gamma_i}$, where $\Gamma_i=\sum_{j=1}^{i}E_j$ and
$(E_j)_{j \geq 1}$ are i.i.d. exponential random variables of mean
$1$. As in example (a) above, we can exclude $0$ from the space.
Clearly, we may assume that the points of $N$ lie in the space
$E=(0,\infty]$.
But the points of $N$ accumulate at
$\infty$, and so, $N$ is {\em not a point process} on $E$: the
number of points in the (relatively compact) set $(x,\infty)$ is
$\infty$, a.s.

(c) The Poisson process on $(0,\infty)$ of intensity $\nu(dx)=\alpha
x^{-\alpha-1}dx$ (for some $\alpha>0$) is included in our framework.
In this case, $N=\sum_{i \geq 1}\delta_{\Gamma_i^{-1/\alpha}}$, and
the points of $N$ accumulate at $0$. This does not contradict the
definition of a point process, since a set of the form
$(0,\varepsilon)$ is {\em not} relatively compact in $(0,\infty)$.
In this case, $N$ is a (well-defined) point process on the space
$E=(0,\infty]$.

}
\end{remark}

Let $M_0=M_p(E) \verb2\2 \{0\}$. For $\mu=\sum_{j \geq 1}
\delta_{t_j} \in M_0$, we let $x_{\mu}:=\sup_{j \geq 1}|t_j|$.
Note that $x_{\mu}<\infty$ for all $\mu \in M_0$. (This is clear
if $E=[-b,a) \cup (a,b]$ for some $0 \leq a<b<\infty$. Suppose now that
$E=[-\infty,-a) \cup (a,\infty]$ for some $0 \leq a<\infty$. Assume that there exists
$\mu=\sum_{j}\delta_{t_j} \in M_0$ with $\sup_j|t_j|=\infty$.
Then, there exists a subsequence $(t_{j_k})_k$ such that $\lim_k
|t_{j_k}|=\infty$. It follows that the point measure $\mu$ has an
infinite number of points in the relatively compact set $[-\infty,c)
\cup (c,\infty]$ for some $c>0$, which is a contradiction.)

Note that for any $\mu \in M_0$, there exists $x \in (0,\infty)$ such that
 $\mu([-x,x]^c)=0$. (We may take $x=x_{\mu}$.)
For any $x>0$, we let
\begin{eqnarray*}
M_x &=& \{\mu \in M_0; \mu([-x,x]^c)>0\}.
\end{eqnarray*}
 Then,
$$\mu \in M_x \quad \mbox{if and only if} \quad x_{\mu}>x.$$

Using (\ref{cond-lambda}) with $B=[-x,x]^c$, the fact that $\mu([-x,x]^c) \geq 1$ for all $\mu \in M_x$, and the function $f(y)=1-e^{-y}$ is non-decreasing, we obtain:
$$\infty>\int (1-e^{-\mu([-x,x]^c)})\lambda(d\mu)=\int_{M_x}(1-e^{-\mu([-x,x]^c)})\lambda(d\mu) \geq (1-e^{-1})\lambda(M_x).$$

Hence $\lambda(M_x)<\infty$ for all $x>0$.
\vspace{3mm}

Recall that $x$ is a {\em fixed atom} of a point process $N$ if
$P(N\{x\}>0)>0$. (Consequently, $x$ is not a fixed atom of $N$, if
$N\{x\}=0$ a.s.) %By Lemma 9.3.III of \cite{daley-verejones08},
It is known that the set $D$ of all fixed atoms of any point process
is countable.

Let $\cM_p(E)$ be the class of Borel sets in $M_p(E)$, and
$\cM_0=\{M \in \cM_p(E);M \subset M_0\}$. We consider the measurable
map $\Phi: M_0 \to (0,\infty)$, defined by $\Phi(\mu)=x_{\mu}.$

\vspace{3mm}

The following theorem is the main result of the present article.

\begin{theorem}
\label{DH-th} Let $(N_{n})_{n \geq 1}$ be a sequence of point
processes on $E$ satisfying (AI). % and (\ref{null-array}).
Suppose
that $N_{i,n} \in M_0$ a.s. for all $n \geq 1$. Let
$Y_{i,n}=\Phi(N_{i,n})$.

Let $N$ be an infinitely divisible point process on $E$ which satisfies
(\ref{cond-N-pm-infty}).
%, and suppose that its canonical measure $\lambda$ satisfies (\ref{cond-lambda-M0}).
Let $D$ be the set of
fixed atoms of $N$ and $D'=\{x>0; x \in D \ \mbox{or} \ -x \in D\}$.
Assume that
\begin{equation}
\label{new-cond-Y} \limsup_{n \to
\infty}\sum_{i=1}^{k_n}P(Y_{i,n}>x)<\infty, \quad \forall x>0, x
\not \in D'.
\end{equation}

Then $N_n \stackrel{d}{\to} N$ if and only if the following two
conditions hold:
\begin{eqnarray*}
(a) & & \sum_{i=1}^{k_n} P(Y_{i,n}>x)  \to \lambda(M_x), \ \mbox{for
any} \
x>0, x \not \in D' \\ %\ \mbox{for which} \ \lambda(M_x)>0\\
(b) & & \sum_{i=1}^{k_n}P(Y_{i,n}>x,N_{i,n} \in M) \to \lambda(M
\cap M_x), \ \mbox{for any} \ x>0, x \not \in
D' \ \mbox{for which} \\
& &  \lambda(M_x)>0, \ \mbox{and for any set} \ M \in \cM_{0} \
\mbox{such that} \ \lambda(\partial M \cap M_x)=0.
\end{eqnarray*}
\end{theorem}

\noindent {\bf Proof:} The ideas of this proof are borrowed from the
proof of Theorem 2.5 of \cite{DH95}. Suppose that $N_n
\stackrel{d}{\to} N$. We first prove (a). Let $\tilde
N_n=\sum_{i=1}^{k_n} \tilde N_{i,n}$ where $(\tilde N_{i,n})_{i \leq
k_n}$ are independent point processes, and $\tilde N_{i,n}
\stackrel{d}{=}N_{i,n}$. Using (AI) and the fact that $N_n
\stackrel{d}{\to} N$, it follows that $\tilde N_n \stackrel{d}{\to}
N$.

Let $x>0,x \not \in D'$ be arbitrary. We prove that (a) holds. By
Lemma 4.4 of \cite{kall83}), $\tilde N_n([-x,x]^c) \stackrel{d}{\to}
N([-x,x]^c)$, since $N\{x\}=N\{-x\}=0$ a.s. Hence,
\begin{equation}
\label{proof-a-step3} -\log P(\tilde N_n([-x,x]^c)=0) \to  -\log
P(N([-x,x]^c)=0).
\end{equation}
Note that % $N([-x,x]^c)$ has a compound Poisson distribution, with
$$E(e^{-uN([-x,x]^c)})=\exp\left\{-\int_{M_0}(1-e^{-u \mu([-x,x]^c)})\lambda(d\mu)
\right\}, \quad \forall u \in \bR_{+},$$ and hence
\begin{equation}
\label{proof-a-step4} -\log P(N([-x,x]^c)=0)= \lambda (\{\mu \in
M_0; \mu([-x,x]^c)>0\})=
\lambda(M_x).%=\nu(x,\infty).
\end{equation}

From (\ref{proof-a-step3}) and (\ref{proof-a-step4}), we obtain
that:
\begin{equation}
\label{proof-a-step1} P(\tilde N_n([-x,x]^c)=0) \to
e^{-\lambda(M_x)}.
\end{equation}

 Let  $a_{i,n}:=P(Y_{i,n}>x)=
P(N_{i,n}([-x,x]^c)>0)$. Using the independence of $(\tilde
N_{i,n})_{i \leq k_n}$, and the fact that $N_{i,n}
\stackrel{d}{=}\tilde N_{i,n}$, we obtain:
\begin{eqnarray}
\nonumber P(\tilde N_n([-x,x]^c)=0)&=&\prod_{i=1}^{k_n} P(\tilde
N_{i,n}([-x,x]^c)=0) =\prod_{i=1}^{k_n} P( N_{i,n}([-x,x]^c)=0) \\
\label{proof-a-step2} &=& \prod_{i=1}^{k_n}(1-a_{i,n}).
\end{eqnarray}

From (\ref{proof-a-step1}) and (\ref{proof-a-step2}), we infer that
$\prod_{i=1}^{k_n}(1-a_{i,n}) \to e^{-\lambda(M_x)}$, and hence
$\sum_{i=1}^{k_n}- \log(1-a_{i,n}) \to \lambda(M_x)$. By
(\ref{null-array}), $\max_{i \leq k_n}a_{i,n} \to 0$. Using
(\ref{new-cond-Y}) and the fact that $-\log(1-x)=x+O(x^2)$ as $x \to
0$, we infer that $\sum_{i=1}^{k_n}a_{i,n} \to \lambda(M_x)$, i.e.
condition (a) holds.

We now prove (b). Let $x>0, x \not \in D'$ be such that
$\lambda(M_x)>0$. Let $P_{n,x}$ and $P_x$ be probability measures on
$M_0$, defined by:
$$P_{n,x}(M)=\frac{\lambda_n(M \cap M_x)}{\lambda_n(M_x)} \quad
\mbox{and} \quad P_{x}(M)=\frac{\lambda(M \cap
M_x)}{\lambda(M_x)},$$ where $\lambda_n=\sum_{i=1}^{k_n}P \circ
N_{i,n}^{-1}$. By (a),
$$\lambda_n(M_x)=\sum_{i=1}^{k_n}P(N_{i,n} \in
M_x)=\sum_{i=1}^{k_n}P(Y_{i,n}>x) \to \lambda(M_x).$$

As in the proof of Theorem 2.5 of \cite{DH95}, using
(\ref{conv-int-lambda-n}), one can prove that: $P_{n,x}
\stackrel{w}{\to} P_{x}$. This is equivalent to saying that
\begin{equation}
\label{conv-M-Pnx} P_{n,x}(M) \to P_{x}(M),
\end{equation} for any $M
\in \cM_0$ with $\lambda(\partial M \cap M_x)=0$. Condition (b)
follows, since
$$P_{n,x}(M)=\frac{\sum_{i=1}^{k_n}P(N_{i,n} \in M \cap M_x)}
{\sum_{i=1}^{k_n}P(N_{i,n} \in
M_x)}=\frac{\sum_{i=1}^{k_n}P(Y_{i,n}>x,N_{i,n} \in
M)}{\sum_{i=1}^{k_n}P(Y_{i,n} >x)}.$$

Suppose that (a) and (b) hold. Using Theorem \ref{kall-th}, it
suffices to prove that (\ref{conv-int-lambda-n}) holds. Let $f \in
C_{K}^{+}(E)$ be arbitrary. Pick $x>0, x \not \in D'$ such that
$\lambda(M_x)>0$ and the support of $f$ is contained in
$[-x,x]^{c}$. (This choice is possible due to
(\ref{compactness-cond2}).)

Let $P_{n,x}$ and $P_{x}$ be defined as above. Let $M \in \cM_0$ be
 such that $P_x(\partial M)=0$. Due to (a) and
(b), (\ref{conv-M-Pnx}) holds. Hence $P_{n,x} \stackrel{w}{\to}
P_{x}$. %By 15.7.6 of \cite{kall83} (since $P_{n,x}$ and $P_{x}$ are
%probability measures), it follows that $P_{n,x} \stackrel{w}{\to}
%P_{x}$.
It follows that for any bounded continuous function $h:M_0 \to \bR$
$$\frac{1}{\lambda_n(M_x)}\int_{M_x}h(\mu)\lambda_n(d\mu) \to
\frac{1}{\lambda(M_x)}\int_{M_x}h(\mu)\lambda(d\mu).$$ Taking
$h(\mu)=e^{-\mu(f)}$ and using the fact that $\lambda_n(M_x) \to
\lambda(M_x)$ (which is condition (a)), we obtain that:
$$\int_{M_x}(1-e^{-\mu(f)})\lambda_n(d\mu) \to
\int_{M_x}(1-e^{-\mu(f)})\lambda(d\mu).$$

\noindent  Note that if $\mu \not \in M_x$, then $\mu([-x,x]^c)=0$,
$\mu(f)=0$ (since the support of $f$ is contained in $[-x,x]^c$),
and hence
$$\int_{M_0}(1-e^{-\mu(f)})\lambda_n(d\mu)= \int_{M_x}(1-e^{-\mu(f)})
\lambda_n(d\mu) $$ $$ \int_{M_0}(1-e^{-\mu(f)})\lambda(d\mu)=
\int_{M_x}(1-e^{-\mu(f)})\lambda(d\mu).$$ Relation
(\ref{conv-int-lambda-n}) follows using the fact that
$\lambda_n(M_0^c)=\lambda(M_0^c)=0$. $\Box$

\begin{remark}
{\rm As it is seen from the proof, in condition (b), one may replace
$M$ by $M'=M \cap M_x$, where $M \in \cM_{0}$ and $ \lambda(\partial
M \cap M_x)=0.$}
\end{remark}

The following result illustrates a typical application of the
Theorem \ref{DH-th}.

\begin{proposition}
\label{DH-th-cor} For each $n \geq 1$, let $(X_{j,n})_{1 \leq j \leq
n}$ be a strictly stationary sequence of random variables with
values in $E$, such that
$$\limsup_{n \to \infty}nP(|X_{1,n}|>\varepsilon)<\infty \quad \mbox{for
all} \
 \varepsilon>0.$$

Suppose that there exists a sequence $(r_n)_n$ of positive integers
such that %$r_n \to \infty$ and
$k_n=[n/r_n] \to \infty$, such that
\begin{equation}\label{aii}
\lim_{n \to \infty}|E(e^{-\sum_{j=1}^{n}f(X_{j,n})})-
\{E(e^{-\sum_{j=1}^{r_n}f(X_{j,n})})\}^{k_n}|=0, \quad \forall f \in
C_{K}^{+}(E).
\end{equation}
 Let $N_{n}=\sum_{j=1}^{n}\delta_{X_{j,n}}$ and
$N_{r_n,n}=\sum_{j=1}^{r_n}\delta_{X_{j,n}}$.

Let $N$ be an infinitely divisible point process on $E$ which satisfies
(\ref{cond-N-pm-infty}).
%, and suppose that its canonical measure $\lambda$ satisfies (\ref{cond-lambda-M0}).
Let $D$ be the set of
fixed atoms of $N$ and $D'=\{x>0; x \in D \ \mbox{or} \ -x
 \in D\}$.

Then $N_{n} \stackrel{d}{\to} N$ if
 and only if the following two conditions hold:
 \begin{eqnarray*}
(a) & & k_n P(\max_{j \leq r_n}|X_{j,n}|>x)  \to \lambda(M_x), \
\mbox{for any} \
x>0, x \not \in D' \\ % \ \mbox{for which} \ \lambda(M_x)>0 \\
(b) & & k_nP(\max_{j \leq r_n}|X_{j,n}|>x,N_{r_n,n} \in M) \to
\lambda(M \cap M_x), \ \mbox{for any} \ x>0, x \not \in
D' \\
& & \mbox{for which} \ \lambda(M_x)>0, \ \mbox{and for any set} \ M
\in \cM_{0} \ \mbox{such that} \ \lambda(\partial M \cap M_x)=0.
\end{eqnarray*}
\end{proposition}

\begin{remark}
{\rm For each $n \geq 1$, let $(X_{j,n})_{1 \leq j \leq n}$ be a
sequence of random variables with values in $E$ fulfilling all the
requirements of Proposition \ref{DH-th-cor}. We claim that if (a)
holds with $r_n=1$ and $k_n=n$ then the limit point process $N$ is a
Poisson process. To see this, let $\mu$ be a measure on $E$ such
that $\lambda(M)=\mu(\{y \in E; \delta_y \in M\}).$ Condition (a) of
Proposition \ref{DH-th-cor} becomes $n P(|X_{1,n}| >x) \to
\mu([-x,x]^c)$ for all $x>0, x \not \in D'$, which is equivalent to
\begin{equation}
\label{resnick-vague} nP(X_{1,n}\in \cdot) \stackrel{v}{\to}
\mu(\cdot).
\end{equation}
Let $N^*_n=\sum_{j=1}^n\delta_{X^*_{j,n}}$ where $(X^*_{j,n})_{j
\leq n}$ are i.i.d. copies of $X_{1,n}$.
 By Proposition 3.21 of \cite{resnick87}, $N^*_n
\stackrel{d}{\to} N$ where $N$ is a Poisson process on $E$ of
intensity $\mu$, such that $N(\{\pm \infty\})=0$ a.s. Combining thi
with (\ref{aii}), we infer that $N_n \stackrel{d}{\to} N$.}
\end{remark}

The next result gives an expression for the limits which appear in
conditions (a) and (b) of Theorem \ref{DH-th}, using a pair
$(\nu,K)$, where $\nu$ is a measure on $(0,\infty)$ and $K$ is a
transition kernel from $(0,\infty)$ to $M_p(E)$ (i.e. $K(x, \cdot)$
is a probability measure on $M_p(E)$ for all $x>0$, and $K(\cdot,M)$
is Borel measurable for any $M \in \cM_p(E)$). In Section
\ref{appl-section}, we will identify the pair $(\nu,K)$ in some
particular cases.

\begin{lemma}
\label{disintegration} Let $\lambda$ be the canonical measure of an
infinitely divisible point process on $E$.
%such that (\ref{cond-lambda-M0}) holds.
Then there exists a measure $\nu$ on
$(0,\infty)$ and a transition kernel $K$ from $(0,\infty)$ to
$M_p(E)$, such that for any $M \in \cM_p(E)$ and for any $x>0$,
\begin{equation}
\label{lambda-nu-K} \lambda(M)=\int_{0}^{\infty}K(y,M)\nu(dy),
\end{equation}
\begin{equation}
\label{rel-nu-K-ab} \lambda(M_x)=\nu(x,\infty) \quad \mbox{and}
\quad \lambda(M \cap M_x)=\int_x^{\infty}K(y,M)\nu(dy).
\end{equation}
\end{lemma}

\noindent {\bf Proof:} Let $\Psi: M_0 \to (0,\infty) \times M_0$ be
the measurable map defined by:
$$\Psi(\mu)=(\Phi(\mu),\mu)=(x_{\mu}, \mu).$$

Let $\Lambda=\lambda \circ \Psi^{-1}$. Since $\lambda$ is
$\sigma$-finite, the measure $\Lambda$ is $\sigma$-finite on
$(0,\infty) \times M_p(E)$. The space $M_p(E)$ is Polish (by 15.7.7
of \cite{kall83}), and hence it has the disintegration property (see
e.g. 15.3.3 of \cite{kall83}). More precisely, for any Borel set $B
\subset (0,\infty)$ and for any set $M \in \cM_p(E)$,
$$\Lambda(B \times M)=\int_{B}K(y,M)\nu(dy),$$ where $\nu=\lambda \circ
\Phi^{-1}$ is the first marginal of $\Lambda$, and $K$ is a
transition kernel from $(0,\infty)$ to $M_p(E)$. In particular,
(\ref{lambda-nu-K}) holds. We have:
\begin{eqnarray*}
\lambda(M_x)&=&\lambda(\{\mu; x_{\mu}>x\})=(\lambda \circ \Phi^{-1})(x,\infty)=\nu(x,\infty)\\
\lambda(M \cap M_x)&=& \lambda(\mu; x_{\mu}>x, \mu \in M)=(\lambda
\circ \Psi^{-1})((x, \infty) \times M)\\
& =& \Lambda((x,\infty) \times M)=\int_{x}^{\infty} K(y,M)\nu(dy).\,\,\,\, \Box
\end{eqnarray*}
The next result is an application of Proposition \ref{DH-th-cor} to
triangular arrays of the form $X_{j,n}=\xi_j/a_n$, where $(\xi_j)_{j
\geq 1}$ is a stationary sequence (combined with Proposition 0.4.(ii) of \cite{resnick87}). This result gives the
relationship between the canonical measure $\lambda$ of the limit
process $N$, and the extremal index of the sequence $(|\xi_j|)_{j
\geq 1}$, when it exists.
\begin{corollary}
\label{extremal-index-cor} Let $(\xi_j)_{j \geq 1}$ be a stationary
sequence of real-valued random variables. Suppose that the extremal
index $\theta$ of the sequence $(|\xi_j|)_{j \geq 1}$ exists and is positive. Let
$(a_n)_{n \geq 1}$ be a sequence of positive real numbers such that
$nP(|\xi_1|> x a_n) \to \mu(x,\infty)$, where $\mu$ is a measure on
$(0,\infty)$, such that $\mu \not = o$ and $\mu(x,\infty)<\infty$ for all $x>0$. Suppose that the triangular array $X_{j,n}=\xi_j/a_n$
satisfies condition (\ref{aii}). Let
$N_n=\sum_{j=1}^{n}\delta_{\xi_j/a_n}$ and
$N_{r_n,n}=\sum_{j=1}^{r_n}\delta_{\xi_j/a_n}$. Let $N$ and $D$ be
as in Proposition \ref{DH-th-cor}.

 Then $N_{n} \stackrel{d}{\to} N$ if
 and only if the following two conditions hold:
 \begin{eqnarray*}
(a) & & k_n P(\max_{j \leq r_n}|\xi_{j}|>xa_n)  \to \theta
\mu(x,\infty), \ \mbox{for any} \
x>0, x \not \in D' \\
(b) & & k_nP(\max_{j \leq r_n}|\xi_{j}|>xa_n,N_{r_n,n} \in M) \to
\lambda(M \cap M_x), \ \mbox{for any} \ x>0 \ \mbox{with} \\
& & \mu(x,\infty)>0, \ \mbox{and for any set} \ M \in \cM_{0} \
\mbox{such that} \ \lambda(\partial M \cap M_x)=0.
\end{eqnarray*}
In this case, $(|\xi_j|)_j$ is regularly varying and $\mu(x,\infty)=x^{-\alpha}$ for some $\alpha>0$.
\end{corollary}

\begin{remark}
{\rm Under the hypothesis of Corollary \ref{extremal-index-cor}, by
Lemma \ref{disintegration}, there exists a transition probability
$K$ from $(0,\infty)$ to $M_p(E)$ such that $\lambda(M)=\theta
\int_0^{\infty} K(y,M)\mu(dy)$. Condition (b) of Corollary
\ref{extremal-index-cor} becomes:
$$k_n\int_x^{\infty}K_n(y,M)P(\max_{j \leq r_n}|\xi_{j}|/a_n\in dy) \to
\theta\int_x^{\infty}K(y,\cdot)\mu(dy),$$ where
$K_n(y,M)=P(N_{r_n,n} \in M| \max_{j \leq r_n}|\xi_{j}|/a_n=y)$. }
%and this explains where the different component of $N$ come from.}
\end{remark}

\section{Applications}
\label{appl-section}

\subsection{Exceedance processes}

In this subsection, we take $E=(0,1]$. Let $(\xi_j)_{j \geq 1}$ be
a stationary sequence of random variables, and $(u_n)_{n}$ be a
sequence of real numbers such that:
\begin{equation}
\label{def-u-n} \lim_{n \to \infty}nP(\xi_1>u_n)=1.
\end{equation}

Supposing that $\xi_j$ represents a measurement made at time $j$, we
define the process $N_n$ which counts the ``normalized'' times $j/n$
when the measurement $\xi_j$ exceeds the level $u_n$, i.e.
\begin{equation}
\label{def-exceedance}
N_n(\cdot)=\sum_{j=1}^{n}\delta_{j/n}(\cdot)1_{\{\xi_j>u_n\}}.
\end{equation} $N_n$ is a point process on
$(0,1]$, called the (time normalized) {\em exceedance process}.

The following mixing-type condition was introduced in \cite{HHL88}. { We
say that $(\xi_j)_{j \geq 1}$ satisfies {\bf condition
$\Delta(\{u_n\})$} if there exists a sequence $(m_n)_{n \geq 1} \subset \bZ_{+}$ such that: $$m_n=o(n) \quad \mbox{and} \quad \alpha_n(m_n) \to 0,$$
where $$\alpha_n(m):=\sup\{|P(A \cap B)-P(A)P(B)|; A \in
\cF_1^k(u_n), \cF_{k+m}^n(u_n), k+m \leq n \},$$ and
$\cF_{i}^j(u_n)=\sigma(\{\{\xi_s \leq u_n\}; i \leq s \leq j\})$.}
Lemma 2.2 of \cite{HHL88} shows that if $(\xi_j)_j$ satisfies
$\Delta(\{u_n\})$, then $(N_n)_{n \geq 1}$ satisfies (AI) with
\begin{equation}
\label{def-N-in} N_{i,n}(\cdot)=\sum_{j \in
J_{i,n}}\delta_{j/n}(\cdot)1_{\{ \xi_j>u_n\}},
\end{equation}
where $J_{i,n}=((i-1)r_n, ir_n]$, $r_n=[n/k_n]$ and $k_n \to
\infty$.

Theorem 4.1 and Theorem 4.2 of \cite{HHL88} show that under
$\Delta(\{u_n\})$, $N_n \stackrel{d}{\to} N$ if and only if the
following two conditions hold:
\begin{eqnarray}
\label{exceedance-cond1} & & k_n P(\max_{j \in J_{1,n}}\xi_j > u_n)
\to  a, \ \mbox{for
some} \ a>0 \\
\label{exceedance-cond2} & & \pi_n(k):=P(\sum_{j \in J_{1,n}}
1_{\{\xi_j>u_n\}}=k|\max_{j \in J_{1,n}} \xi_j>u_n) \to \pi_k, \
\forall  k \geq 1.
\end{eqnarray}

\noindent In this case, %$a$ is the extremal index of $(\xi_j)_j$ and
$N$ is a compound Poisson process on $(0,1]$ with Poisson rate $a$
and the distribution of multiplicities $(\pi_k)_{k \geq 1}$, i.e.
$$L_{N}(f)=\exp \left\{-a\int_{0}^{1}(1-\sum_{k \geq 1}
e^{-kf(x)})\pi_k dx \right\}.$$

Note that condition
$\Delta(\{u_n\})$ is a mixing-type condition which is sufficient
for (AI) (by Lemma 2.2 in [4]), but may not be necessary.
However, as our next result shows, condition $\Delta(\{u_n\})$ is
much stronger than needed for the convergence of $(N_n)_{n \geq 1}$.
In fact, this convergence can be obtained under (AI) alone.

\begin{proposition}
\label{exceed-AI-prop}
Suppose that $(N_n)_{n \geq 1}$ satisfies (AI) with $(N_{i,n})_{i
\leq k_n}$ given by (\ref{def-N-in}). If (\ref{exceedance-cond1})
and (\ref{exceedance-cond2}) hold, then conditions (a) and (b) of
Theorem \ref{DH-th} are satisfied, and (\ref{rel-nu-K-ab}) holds
with
$$\nu(dy)=a 1_{(0,1]}(y)dy \quad \mbox{and} \quad K(y,M)=\sum_{k \geq 1}
\pi_k 1_{M}(k \delta_y).$$
\end{proposition}

Before proving Proposition \ref{exceed-AI-prop}, we shall discuss an application of Proposition \ref{exceed-AI-prop} to the case of associated random variables $(\xi_j)_{j \geq 1}$, for which the dependence condition $\Delta(\{u_n\})$ is not appropriate.
Recall that the random variables $(\xi_{j})_{j \geq 1}$ are {\em
associated} if
$${\rm Cov}(g(\xi_{1},\cdots, \xi_n), h(\xi_{1},\cdots,\xi_n)) \geq 0,$$
for any $n \geq 1$, and any coordinate-wise non-decreasing functions
$g:\bR^{n} \to \bR$ and $h: \bR^{n} \to \bR$ for which the covariance is well-defined (see e.g. \cite{barlow-prochan81}, \cite{EPW67}). This notion is very different from mixing. To see this, let $(\varepsilon_i)_{i \in \bZ}$ be a sequence of i.i.d. Bernoulli random variables with parameter 1/2. The linear process $$\xi_j=\sum_{i \geq 0} 2^{-i}\varepsilon_{j-i}, \quad j \in \bZ$$ is associated (by $\cP_2$ and $\cP_4$ of \cite{EPW67}), but fails to be mixing (see \cite{bradley86} and the references therein).

The next result identifies a condition under which the sequence $(N_n)_{n \geq 1}$ of exceedance processes satisfies (AI), when $(\xi_j)_{j \geq 1}$ are associated.

\begin{lemma}\label{asso}
Let $(\xi_j)_{j\geq 1}$ be a stationary sequence of associated random variables, and $(u_n)_{n \geq 1}$ be a sequence of real numbers such that (\ref{def-u-n}) holds. Suppose that
\begin{equation}
\label{assoc-cond}
\lim_{m \rightarrow \infty}\limsup_{n \rightarrow \infty}n\sum_{j=m+1}^n {\rm
Cov}(1_{\{\xi_1>u_n\}}, 1_{\{\xi_j>u_n\}})=0.
\end{equation}

Then the sequence $(N_n)_{n \geq 1}$ of exceedence processes defined by (\ref{def-exceedance}) satisfies (AI) (with $(N_{i,n})_{i \leq k_n}$ defined by (\ref{def-N-in})), assuming that $\lim_{n \rightarrow \infty}r_n=\infty$.
\end{lemma}

\noindent {\bf Proof of Lemma \ref{asso}:} Let $f \in C_{K}^{+}(E)$ be arbitrary. Without loss of generality, we may assume that $f(x) \leq 1$ for all $x \in E$. For each $n \geq 1$, define
$$Y_{j,n}=f(j/n)1_{\{ \xi_j >u_n\}}, \quad 1 \leq j \leq n.$$
Note that the random variables $(Y_{j,n})_{j \leq n}$ are associated. Clearly,
$$|E(e^{-N_n(f)})-
\prod_{i=1}^{k_n}E(e^{-N_{i,n}(f)})|= |E(e^{-\sum_{j=1}^nY_{j,n}})-
\prod_{i=1}^{k_n}E(e^{-\sum_{j\in J_{i,n} }Y_{j,n}})|.$$

We follow the lines of proof of Lemma 5.4 of \cite{BL08}.
For each $1 \leq i \leq k_n$, let $H_{i,n}^{(m)}$ be the (big) block
of consecutive integers between $(i-1)r_n+1$ and $ir_n-m$ and
$I_{i,n}^{(m)}$ be the (small) block of size $m$, consisting of
consecutive integers between $ir_n-m$ and $ir_n$. Let
$$U_{i,n}^{(m)}=\sum_{j \in H_{i,n}^{(m)}}Y_{j,n}.$$

Similarly to (32) and (33) of \cite{BL08}, in order to prove (AI),
it suffices to show that
\begin{equation}
\label{AD-third-conv-2} \lim_{m \to \infty} \limsup_{n \to \infty}
|E(e^{-\sum_{i=1}^{k_n}U_{i,n}^{(m)}})-\prod_{i=1}^{k_n} E(e^{-U_{i,n}^{(m)}})|
=0.
\end{equation}
For this we argue as for (36) and (37) in \cite{BL08} to get,
\begin{eqnarray*}
\label{association-induction}
&& |E(e^{-\sum_{i=1}^{k_n}U_{i,n}^{(m)}})-\prod_{i=1}^{k_n} E(e^{-U_{i,n}^{(m)}})|
\leq \sum_{1 \leq i <l \leq k_n}\sum_{j \in H_{i,n}^{(m)}}
\sum_{j' \in H_{l,n}^{(m)}} {\rm Cov}(Y_{j,n},Y_{j',n})\\
&& =\sum_{1 \leq i <l \leq k_n}\sum_{j \in H_{i,n}^{(m)}}
\sum_{j' \in H_{l,n}^{(m)}}f(j/n)f(j'/n){\rm
Cov}(1_{\xi_j>u_n},1_{\xi_{j'}>u_n})\\
&& \leq 2n\sum_{l=m+1}^n{\rm Cov}(1_{\xi_1>u_n},1_{\xi_{l}>u_n}),
\end{eqnarray*}
and relation (\ref{AD-third-conv-2}) follows from (\ref{assoc-cond}). $\Box$
\\
\\
\noindent {\bf Proof of Proposition \ref{exceed-AI-prop}:} In this case, $D'=\emptyset$ (since $N$ does
not have fixed atoms), and
%The canonical measure $\lambda$ of $N$ is given by
%(\ref{def-lambda-cPoisson}), with $\pi_k(x)=\pi_k$ and $\nu^*(dy)=a
%dy$, i.e.
%$$\lambda(M)=a\sum_{k \geq 1} \pi_k  \int_0^1 1_{M}(k \delta_y)
%dy=\int_{0}^{1}K^*(y,M)\nu^*(dy).$$ Clearly, $\lambda(M_0^c)=0$.
$$Y_{i,n}=\max\{j/n; j \in J_{i,n}, \xi_j>u_n \}, \quad
\mbox{for all} \quad i=1, \ldots,k_n.$$ To verify condition (a) of
Theorem \ref{DH-th}, we show that for any $x \in (0,1]$:
\begin{equation}
\label{verif-cond-a} \sum_{i=1}^{k_n}P(Y_{i,n}>x) \to
a(1-x)=\nu(x,\infty).
\end{equation}

Let $x \in (0,1]$ be arbitrary. Then $x \leq k_n r_n/n$ for $n$
large enough (since $k_n r_n/n \to 1$). It follows that for $n$
large enough, $nx \in (0,k_n r_n] =\cup_{i=1}^{k_n}J_{i,n}$, and
there exists $i_n \leq k_n$ such that $nx \in J_{i_n,n}$, i.e.
$(i_n-1)r_n/n+1/n \leq x<i_nr_n/n$. We write
\begin{eqnarray*}
\sum_{i=1}^{k_n}P(Y_{i,n}>x)&=&\sum_{i=1}^{i_n-1}P(Y_{i,n}>x)
+P(Y_{i_n,n}>x)+\sum_{i=i_n+1}^{k_n}P(Y_{i,n}>x)\\
&=:&I_1(n)+I_2(n)+I_3(n).
\end{eqnarray*}

Note that $\{Y_{i,n} >x\}=\emptyset$ for all $i \leq i_{n}-1$, and
hence $ I_1(n)=0$.
%\begin{equation}
%\label{verif-cond-a-step1} I_1(n)=0.
%\end{equation}

%(To see this, assume that $Y_{i,n}>x$ and $i \leq i_n-1$. Then there
%exists $j \in J_{i,n}$ such that $\xi_j>u_n$ and $j/n>x$. On one
%hand, we have $x<j/n<ir_n/n \leq (i_n-1)r_n/n$, and on the other
%hand $x >(i_n -1)r_n/n$, since $nx \in J_{i_n,n}$. This is a
%contradiction.)

We have $\{Y_{i_n,n}>x\} \subset \cup_{j \in
J_{i_n,n}}\{\xi_j>u_n\}$ and hence,  by (\ref{def-u-n}), $I_2(n)
\leq \sum_{j \in J_{i_n,n}} P(\xi_j>u_n)=r_n P(\xi_1>u_n) \leq
\frac{1}{k_n} \cdot n P(\xi_1>u_n) \to 0$.
%\begin{equation}
%\label{verif-cond-a-step2} I_2(n) \leq \sum_{j \in J_{i_n,n}}
%P(\xi_j>u_n)=r_n P(\xi_1>u_n) \leq \frac{1}{k_n} \cdot n
%P(\xi_1>u_n) \to 0.
%\end{equation}

Finally, we claim that $\{Y_{i,n} >x\}= \{\max_{j \in J_{i,n}}
\xi_j>u_n\}$ for any $i \geq i_{n}+1$. This is true since for any $j
\in J_{i,n}$ and $i \geq i_n+1$, we have:
$$\frac{j}{n} \geq \frac{(i-1)r_n+1}{n} \geq \frac{i_n
r_n}{n}+\frac{1}{n}>x+\frac{1}{n}>x.$$

\noindent Using the stationary of the sequence $(\xi_j)_j$, the fact
that $i_n/k_n \to x$, and (\ref{exceedance-cond1}), we obtain:
%\begin{equation}
%\label{verif-cond-a-step3}
$$I_3(n)= \sum_{i=i_n+1}^{k_n}P(\max_{j
\in J_{i,n}}\xi_j>u_n)=\left(1-\frac{i_n}{k_n}\right)  k_n P(\max_{j
\in J_{1,n}} \xi_j>u_n) \to (1-x)a.$$
%\end{equation}
Relation (\ref{verif-cond-a}) follows
%from
%(\ref{verif-cond-a-step1}), (\ref{verif-cond-a-step2}) and
%(\ref{verif-cond-a-step3}).
%
To verify condition (b) of Theorem \ref{DH-th}, we will show that
for any $x \in (0,1]$ and for any $M \in \cM_0$,
$$\sum_{i=1}^{k_n}P(N_{i,n} \in M, Y_{i,n}>x) \to
 a\sum_{k \geq 1} \pi_k \int_x^1 1_{M}(k \delta_y)dy=\int_x^{\infty}
K(y,M)\nu(dy).$$

Let $x \in (0,1]$ and $M \in \cM_0$ be arbitrary. Then $nx \in
J_{i_n,n}$ for some $i_n \leq k_n$, and $n$ large enough. We write
$\sum_{i=1}^{k_n}P(N_{i,n} \in M, Y_{i,n}>x)=J_1(n)+J_2(n)+J_3(n)$,
where
\begin{eqnarray*}
J_1(n)&:=& \sum_{i=1}^{i_n-1}P(N_{i,n} \in M, Y_{i,n}>x) \leq
I_1(n)=0 \\
J_2(n)&:=& P(N_{i_n,n} \in M, Y_{i_n,n}>x) \leq I_2(n) \to 0 \\
J_3(n)&:=&\sum_{i=i_n+1}^{k_n}P(N_{i,n} \in M,
Y_{i,n}>x)=\sum_{i=i_n+1}^{k_n}P(N_{i,n} \in M, \max_{j \in
J_{i,n}}\xi_j>u_n)
\end{eqnarray*}

Hence, it suffices to show that:
\begin{equation}
\label{verif-cond-b} J_3'(n):=\sum_{i=i_n+1}^{k_n}P(N_{i,n} \in M,
\max_{j \in J_{i,n}}\xi_j>u_n) \to  a\sum_{k \geq 1} \pi_k \int_x^1
1_{M}(k \delta_y)dy.
\end{equation}

By the stationary of the sequence $(\xi_j)_{j}$
\begin{eqnarray}
\nonumber J_3'(n)&=& \sum_{i=i_n+1}^{k_n} \sum_{k=1}^{r_n}P(N_{i,n}
\in M, \max_{j \in J_{i,n}} \xi_j>u_n, \sum_{j \in J_{i,n}}
1_{\{\xi_j>u_n\}}=k)
\\
\nonumber &=&  \sum_{k=1}^{r_n}\sum_{i=i_n+1}^{k_n}P(N_{i,n} \in
M|\sum_{j \in J_{i,n}} 1_{\{\xi_j>u_n\}}=k) P(\max_{j \in J_{i,n}}
\xi_j>u_n)\pi_n(k) \\
&=& \label{calcul-J3} \frac{n}{r_n}P(\max_{j \leq r_n} \xi_j>u_n)
\sum_{k=1}^{r_n}p_n(k,M)\pi_n(k),
\end{eqnarray}
where
$$p_n(k,M):=\frac{r_n}{n}\sum_{i=i_n+1}^{k_n}P(N_{i,n} \in M|
\sum_{j \in J_{i,n}} 1_{\{\xi_j>u_n \}}=k).$$

Let $\cD$ be the class of sets $M \in \cM_0$, for which there exists
a set $I_M \subset \{1, \ldots,k_n\}$ (which does {\em not} depend
on $n$) such that for any $n \geq 1$, for any $i \in \{1, \ldots,
k_n\} \verb2\2 I_M$ and for any $k \in \{1, \ldots, r_n\}$, we have:
\begin{eqnarray}
\nonumber \lefteqn{ P(\sum_{j \in J_{i,n}} \delta_{j/n}(\cdot)
1_{\{\xi_j>u_n\}} \in M|\sum_{ j \in J_{i,n}} 1_{\{\xi_j>u_n\}}=k)=}
\\
\label{dirac-equality} & & P( \sum_{j \in J_{i,n}} \delta_{y}(\cdot)
1_{\{\xi_j>u_n\}} \in M|\sum_{ j \in J_{i,n}} 1_{\{\xi_j>u_n\}}=k),
\quad \forall y \in J_{i,n}/n.
\end{eqnarray}

Using a monotone class argument, one can prove that:
\begin{equation}
\label{D-equal-M0} \cD=\cM_0.
\end{equation}

This shows that if $i  \in \{1, \ldots,k_n\} \verb2\2 I_M$, then
\begin{eqnarray*}
f_{i,n}(y)&:=&P(\sum_{j \in J_{i,n}} \delta_y(\cdot) 1_{\{
\xi_j>u_n\}} \in M|\sum_{j \in J_{i,n}} 1_{\{\xi_j>u_n \}}=k) \\
&=& P(\sum_{j \in J_{i,n}} \delta_{j/n}(\cdot) 1_{\{
\xi_j>u_n\}} \in M|\sum_{j \in J_{i,n}} 1_{\{\xi_j>u_n \}}=k) \\
&=& P(N_{i,n} \in M|\sum_{j \in J_{i,n}} 1_{\{\xi_j>u_n \}}=k),
\quad \mbox{for all} \ y \in J_{i,n}/n,
\end{eqnarray*}
and hence
$$\int_{J_{i,n}/n}f_{i,n}(y)dy=\frac{r_n}{n}
P(N_{i,n} \in M|\sum_{j \in J_{i,n}} 1_{\{\xi_j>u_n \}}=k).$$ On the
other hand, $f_{i,n}(y)=1_{M}(k \delta_y)$ and hence
$$\int_{J_{i,n}/n}f_{i,n}(y)dy=\int_{J_{i,n}/n}1_{M}(k \delta_y)dy.$$
It follows that if $i \in \{1,\ldots,k_n\} \verb2\2 I_M$,
$$\frac{r_n}{n}
P(N_{i,n} \in M|\sum_{j \in J_{i,n}} 1_{\{\xi_j>u_n
\}}=k)=\int_{J_{i,n}/n}1_{M}(k \delta_y)dy.$$

We return now to the calculation of $p_{n}(k,M)$. We split the sum
over $i=i_n+1, \ldots,k_n$ into two sums, which contain
 the terms corresponding to indices $i \in I_M$, respectively $i
\not \in I_M$. The second sum is bounded by ${\rm card}(I_M)r_n/n$.
More precisely, we have:
\begin{eqnarray}
\nonumber p_n(k,M)&=& \sum_{i \not \in I_M}\frac{r_n}{n} P(N_{i,n}
\in M|\sum_{j \in
J_{i,n}} 1_{\{\xi_j>u_n \}}=k)+O\left(\frac{r_n}{n} \right)\\
&=& \sum_{i \not \in I_M} \int_{J_{i,n}/n} 1_{M}(k
\delta_y)dy +O\left(\frac{r_n}{n} \right) \\
\nonumber &=&  \sum_{i=i_n+1}^{k_n} \int_{J_{i,n}/n} 1_{M}(k
\delta_y)dy-\sum_{i_n+1 \leq i \leq k_n, i \in I_M} \int_{J_{i,n}/n}
1_{M}(k \delta_y)dy +O\left(\frac{r_n}{n} \right) \\
&=& \int_{i_n r_n/n}^{k_n r_n/n} 1_{M}(k \delta_y) dy
+O\left(\frac{r_n}{n} \right)
= \int_{x}^{1} 1_{M}(k
\delta_y)dy+O\left(\frac{r_n}{n}\right).\label{calcul-p}
\end{eqnarray}

Using (\ref{calcul-J3}) and (\ref{calcul-p}), we obtain:
$$J_3'(n)=\frac{n}{r_n}P(\max_{j \leq r_n}\xi_j>u_n)
 \left\{ \sum_{k=1}^{r_n} \pi_n(k)\int_x^1 1_{M}(k \delta_y)dy+
 O\left( \frac{r_n}{n}\right) \right\}$$

Using (\ref{exceedance-cond1}) and (\ref{exceedance-cond2}), it
follows that:
$J_3'(n) \to a \sum_{k \geq 1}\pi_k \int_x^{1}
 1_{M}(k \delta_y)dy.$
Here, we used the fact that $$\left|\sum_{k=1}^{r_n}
\pi_n(k)\int_x^1 1_{M}(k \delta_y)dy-\sum_{k \geq 1}\pi_k \int_x^{1}
 1_{M}(k \delta_y)dy\right| \leq
  (1-x)\sum_{k \geq 1}|\pi_n(k)-\pi_k| \to 0,$$
where the last convergence is justified by Scheff\'e's theorem (see
e.g. Theorem 16.12 of \cite{billingsley95}), since $\sum_{k \geq
1}\pi_n(k)=\sum_{k \geq 1} \pi(k)=1$. This concludes the proof of
(\ref{verif-cond-b}).  $\Box$

\subsection{Processes based on regularly varying sequences}

In this subsection, we take $E=\bar \bR \verb2\2 \{0\}$. Let
$(\xi_j)_{j \geq 1}$ be a stationary sequence of random variables
with values in $\bR \verb2\2 \{0\}$, and $(a_n)_{n}$ be a sequence
of positive real numbers such that:
\begin{equation}
\label{def-a-n} \lim_{n \to \infty}nP(|\xi_1|>a_n)=1.
\end{equation}

We consider the following point process on $\bar \bR \verb2\2
\{0\}$: $N_n=\sum_{j=1}^{n} \delta_{\xi_j/a_n}.$ The mixing
condition introduced in \cite{DH95} is the following. {\rm We say
that $(\xi_j)_{j \geq 1}$ satisfies {\bf condition $\cA(\{a_n\})$}
if there exists a sequence $k_n \to \infty$ with $r_n=[n/k_n] \to
\infty$ such that:
$$\lim_{n \to \infty} |E(e^{-\sum_{j=1}^{n}f(\xi_j/a_n) })-
\{E(e^{-\sum_{j=1}^{r_n}f(\xi_j/a_n)} \}^{k_n} |=0, \quad \forall f
\in C_{K}^{+}(E).$$ }

Note that $\cA(\{a_n\})$ is equivalent to saying that $(N_n)_{n \geq
1}$ satisfies (AI), with
\begin{equation}
\label{def-N-in-RV} N_{i,n}=\sum_{j \in J_{i,n}}\delta_{\xi_j/a_n},
\end{equation}
and $J_{i,n}=((i-1)r_n ,i r_n]$. In particular, if $(\xi_j)_j$ is
strongly mixing, then $\cA(\{a_n\})$ holds (see e.g. Lemma 5.1 of
\cite{BL08}). In addition, suppose that $\xi_1$ has {\em regularly
varying tail probabilities} of order $\alpha >0$, i.e.
\begin{equation}
\label{def-reg-var} P(|\xi_1|>x)=x^{-\alpha}L(x), \quad \lim_{x \to
\infty}\frac{P(\xi_1>x)}{P(|\xi_1|>x)}=p, \quad \lim_{x \to
\infty}\frac{P(\xi_1<-x)}{P(|\xi_1|>x)}=q,
\end{equation} where $L$
is a slowly varying function, $p \in [0,1]$ and $q=1-p$.

Let $\tilde M=\{ \mu \in M_p(E); \mu([-1,1]^c)=0,
\mu(\{-1,1\})>0\}$, and consider the following point process with
values in $\tilde M$,
$\xi_{1,n}=\sum_{j=1}^{r_n}\delta_{\xi_j/\max_{i \leq r_n}|\xi_i|}.$

Theorem 2.3 and Theorem 2.5 of \cite{DH95} show that under
(\ref{def-reg-var}) and $\cA(\{a_n\})$, $N_n \stackrel{d}{\to} N$ if
and only if the following two conditions hold:
\begin{eqnarray}
\label{reg-var-cond1} & & k_n P(\max_{j \leq r_n}|\xi_j|>a_n
x) \to \theta x^{-\alpha}, \quad \forall x>0,  \\
\label{reg-var-cond2} & &  P(\xi_{1,n} \in \cdot \ |\max_{j \leq
r_n}|\xi_j|>a_n x) \stackrel{w}{\to} {\cal Q}(\cdot), \quad \forall
x>0,
\end{eqnarray}
where $\theta \in (0,1]$ and ${\cal Q}$ is a probability measure on
$\tilde M$. In this case, $\theta$ is the extremal index of
$(|\xi_j|)_{j \geq 1}$, and $N$ is an ID point process (without
fixed atoms).
%, whose canonical measure $\lambda$ satisfies (\ref{cond-lambda-M0}).
 Since $\lambda \circ \Omega^{-1}=\nu \times
{\cal Q}$, where
\begin{equation}
\label{def-nu} \nu(dy)=\theta \alpha y^{-\alpha-1}dy
\end{equation}
 and $\Omega:
M_0 \to (0,\infty) \times \tilde M$ is defined by
$\Omega(\mu)=(x_{\mu}, \mu(x_{\mu} \cdot))$, it follows that:
\begin{equation}
\label{Laplace-reg-var} L_{N}(f) =\exp
\left\{-\int_{0}^{\infty}\int_{\tilde M} (1-e^{-\mu( f(y \cdot))})
{\cal Q}(d\mu)\nu(dy) \right\}.
\end{equation}

On the other hand, by Lemma \ref{disintegration}, the Laplace
functional of $N$ is given by:
$$L_{N}(f)=\exp \left\{-\int_{0}^{\infty} \int_{M_0}
(1-e^{-\mu(f)})K(y,d\mu)\nu(dy).
 \right\}$$
This allows us to identify the relationship between the probability
measure ${\cal Q}$ and the kernel $K$, namely:
\begin{equation}
\label{def-K} K(y,M)={\cal Q}(\pi_y^{-1}(M) \cap \tilde M),
\end{equation} where $\pi_y:M_0 \to M_0$ is defined by
$\pi_y(\mu)=\mu(y^{-1} \ \cdot \ )$, i.e.
$\pi_y(\sum_{j}\delta_{t_j})=\sum_{j}\delta_{yt_j}$.

The next result shows that conditions (\ref{reg-var-cond1}) and
(\ref{reg-var-cond2}) are equivalent to conditions (a) and (b) of
Theorem \ref{DH-th}, in the setting of the present subsection.

\begin{proposition}
Suppose that the sequence $(N_n)_n$ satisfies (AI) with
$(N_{i,n})_{i \leq k_n}$ given by (\ref{def-N-in-RV}). Then
 (\ref{reg-var-cond1}) and
(\ref{reg-var-cond2}) are equivalent to conditions (a) and (b) of
Theorem \ref{DH-th}, with $\lambda(M_x),\lambda(M \cap M_x)$ given
by (\ref{rel-nu-K-ab}), and $\nu,K$ given by (\ref{def-nu}) and
(\ref{def-K}).
\end{proposition}

\noindent {\bf Proof:} In this case, $D'=\emptyset$ and
%(since $N$ has no fixed atoms) and
$Y_{1,n}=\max_{j \leq r_n}|\xi_j|/a_n$. Condition (a) of Theorem
\ref{DH-th} is in fact (\ref{reg-var-cond1}). Assuming now that (a)
holds, we show that condition (b) of Theorem \ref{DH-th} is
equivalent to (\ref{reg-var-cond2}).

Let $\lambda_n=k_n P \circ N_{1,n}^{-1}$. Define the following
probability measures on $M_0$:
\begin{eqnarray*}
P_{n,x}(M)&=&\frac{\lambda_n(M \cap
M_x)}{\lambda_n(M_x)}=\frac{k_n P(N_{1,n} \in M,Y_{1,n}>x)}{k_n P(Y_{1,n}>x)} \\
P_{x}(M)&=&\frac{\lambda(M \cap
M_x)}{\lambda(M_x)}=\frac{1}{\nu(x,\infty)}
\int_x^{\infty}K(y,M)\nu(dy).
 \end{eqnarray*}

By (a), $k_n P(Y_{1,n}>x) \to \nu(x,\infty)$. Hence, condition (b)
of Theorem \ref{DH-th} is equivalent to

\begin{equation}
\label{equiv-b-step1} P_{n,x} \stackrel{w}{\to} P_{x}, \quad
\mbox{for all} \ x>0.
\end{equation}

On the other hand, for any $M \in \tilde \cM$,
\begin{eqnarray*}
{\lefteqn{k_n P(\xi_{1,n} \in M,Y_{1,n}>x)= k_n P(N_{1,n} \in \{\mu \in M_0;
x_{\mu}>x, \mu(x_{\mu} \cdot) \in M \}, N_{1,n} \in M_x) }}\\
&&= k_n (P \circ N_{1,n}^{-1})(\{\mu \in M_0; x_{\mu}>x, \mu(x_{\mu}
\cdot) \in M \}\cap M_x)
=\lambda_n( \Omega^{-1}((x,\infty) \times M) \cap M_x),
\end{eqnarray*}
and hence
$$P(\xi_{1,n} \in M|Y_{1,n}>x)=
\frac{\lambda_n( \Omega^{-1}((x,\infty) \times M) \cap
M_x)}{\lambda_n(M_x)}=(P_{n,x} \circ \Omega^{-1})((x,\infty) \times
M).$$

\noindent Since $\lambda \circ \Omega^{-1}=\nu \times {\cal Q}$,
$${\cal Q}(M)=\frac{\lambda (\Omega^{-1}((x,\infty) \times M))}{\nu^*(x,\infty)}
=(P_{x} \circ \Omega^{-1})((x,\infty) \times M), \quad \forall M \in
\tilde M.$$

\noindent Therefore, relation (\ref{reg-var-cond2}) is equivalent
to:
\begin{equation}
\label{equiv-b-step2} P_{n,x} \circ \Omega^{-1} ((x,\infty) \times
\cdot) \stackrel{w}{\to} P_{x} \circ \Omega^{-1} ((x,\infty) \times
\cdot) , \quad \mbox{for all} \ x>0.
\end{equation}

Using the argument on page 888 of \cite{DH95}, one can show that
(\ref{equiv-b-step1}) is equivalent to (\ref{equiv-b-step2}). (This
argument uses the continuous mapping theorem, and the fact that both
$\Omega$ and $\Omega^{-1}$ are continuous.) $\Box$

\end{document}